\documentclass[12pt]{article}
\usepackage{amssymb}

\newtheorem{thm}     {Theorem}[section]
\newtheorem{prop}    [thm]{Proposition}

\newtheorem{cor}     [thm]{Corollary}
\newtheorem{lemma}   [thm]{Lemma}

\newcommand{\proof} {\noindent{\bf Proof. }}

\newcommand{\B}{\mathbb B}
\newcommand{\C}{\mathbb C}
\newcommand{\D}{\mathbb D}

\newcommand{\R}{\mathbb R}

\newcommand{\st}{{\rm st}}

\def\Re{{\rm Re\,}}
\def\Im{{\rm Im\,}}
\def\bar{\overline}

\begin{document}

\title{Pluripolar sets, real submanifolds and pseudoholomorphic discs }
\author{Alexandre Sukhov{*} }
\date{}
\maketitle

{\small

* Universit\'e des Sciences et Technologies de Lille, Laboratoire
Paul Painlev\'e,
U.F.R. de
Math\'e-matique, 59655 Villeneuve d'Ascq, Cedex, France, sukhov@math.univ-lille1.fr
The author is partially suported by Labex CEMPI.

Instiut of Mathematics with Computing Centre - Subdivision of the Ufa Research Centre of Russian
Academy of Sciences, 45008, Chernyshevsky Str. 112, Ufa, Russia.

}
\bigskip

{\small Abstract. We prove that a compact subset of full measure on a generic submanifold of an almost complex manifold is not a pluripolar set. Several related results on boundary behavior of plurisubharmonic functions are established. Our approach is based on gluing a family of complex discs to a generic manifold along a boundary arc and could admit further applications.}



\bigskip

\section{Introduction}

The foundations of the theory of almost complex structures go back to the classical  work of Newlander-Nirenberg (see for example, \cite{Aud}) where a complete criterion of integrability of these structures is established. The modern period of development has began after the famous paper by Gromov \cite{Aud} who discovered a deep connection between the almost complex and symplectic geometry. Since the analysis on almost complex manifolds became a powerful tool in the symplectic geometry.

From an analytic point of view (which is in the focus  of the present paper), the analysis on almost complex manifolds has several features unusual in the  much better understood case of integrable almost complex structures. One of them is that for a "generic" choice of an almost complex structure of complex dimension $\ge 2$, the only holomorphic functions are constants (even locally). By contrary, there are two objects (they can be viewed as dual ones) which always exist, at least, locally: pseudoholomorphic discs and plurisubharmonic functions. Both of them   are important technical tools in the symplectic geometry. Theory of pseudoholomorphic discs is (relatively) well elaborated now (although there are many interesting open questions remaining). Theory of plurisubharmonic functions on almost complex manifolds is much younger and its development is rather recent; many quite natural results remain open.

The goal of the present paper is to study some boundary properties of plurisubharmonic functions near 
real submanifolds of almost complex manifolds. Our main inspiration is the well-known paper by A.Sadullaev \cite{Sa1} where he established several useful results on boundary  uniqueness for plurisubharmonic functions as well as the two constant type theorems in $\C^n$. His main technical tool is a construction (due to S.Pinchuk \cite{Pi}) of local family of holomorphic discs glued along the upper semi-circle to a prescribed generic totally real manifold  in $\C^n$. In the present paper we extend these results to the almost complex case. This first step was done recently in \cite{Su1} but here we present a considerably more advanced results. We hope that the almost complex analog of the Pinchuk-Sadullaev gluing discs construction elaborated in the present paper  will have other applications.

The organization of the  paper can be seen from the contents. Sections 2 and 3 are preliminary. Section 4 contains the main technical tool (the construction of pseudoholomorphic discs in the spirit of \cite{Pi,Sa1}). The main results are contained in Section 5.

\tableofcontents

\section{ Almost complex manifolds  and pseduholomorphic discs} In this section we recall basic notions of the almost complex geometry making our presentation more convenient for specialists in Analysis. Everywhere through this paper we assume that manifolds and almost complex structures are of class $C^\infty$; notice the main results remain true under considerably weaker regularity assumptions.

\subsection{ Almost complex manifolds} Let $M$ be a smooth  manifold of dimension $2n$. {\it An almost complex structure} $J$ on $M$ is a smooth map  which associates to every point $p \in M$ a linear isomorphism $J(p): T_pM \to T_pM$ of the tangent space $T_pM$ such that $J(p)^2 = -I$; here  $I$ denotes the identity map of $T_pM$. Thus, every $J(p)$ is a linear complex structure on $T_pM$. A couple $(M,J)$ is called {\it an almost complex manifold} of complex dimension n. Note that every almost complex manifold admits the canonical orientation represented by $(e_1,Je_1,....,e_n,Je_n)$ where $(e_1,....,e_n)$ is any complex basis of $(T_pM,J(p))$.

An important example is provided by the {\it standard complex structure} $J_{st} = J_{st}^{(2)}$ on $M = \R^2$  which is given in the canonical coordinates of $\R^2$ by the matrix 

\begin{eqnarray}
\label{J_st}
J_{st}= \left(
\begin{array}{cll}
0 & & -1\\
1 & & 0
\end{array}
\right)
\end{eqnarray}
More generally, the standard complex structure $J_{st}$ on $\R^{2n}$ is represented by the block diagonal matrix $diag(J_{st}^{(2)},...,J_{st}^{(2)})$ (usually we drop the notation of dimension because its value  will be clear from the context). As usual,  setting $iv := Jv$ for $v \in \R^{2n}$, we identify $(\R^{2n},J_{st})$ with $\C^n$ using the notation 
$z = x + iy = x + Jy$ for the standard complex coordinates $z = (z_1,...,z_n) \in \C^n$.

Let $(M,J)$ and $(M',J')$ be smooth  almost complex manifolds. A $C^1$-map $f:M' \to M$ is called  
$(J',J)$-complex or  $(J',J)$-holomorphic  if it satisfies {\it the Cauchy-Riemann equations} 
\begin{eqnarray}
\label{CRglobal}
df \circ J' = J \circ df.
\end{eqnarray}

Note that   a map $f: \C^n \to \C^m$ is $(J_{st},J_{st})$-holomorphic if and only if each component of $f$ is a usual holomorphic function.

 Every almost complex manifold
$(M,J)$ can be viewed locally as the unit ball $\B$ in
$\C^n$ equipped with a small almost complex
deformation of $J_{st}$. The following statement is often very useful.
\begin{lemma}
\label{lemma1}
Let $(M,J)$ be an almost complex manifold. Then for every point $p \in
M$, every  $m \geq 0$ and   $\lambda_0 > 0$ there exist a neighborhood $U$ of $p$ and a
coordinate diffeomorphism $z: U \rightarrow \B$ such that
$z(p) = 0$, $dz(p) \circ J(p) \circ dz^{-1}(0) = J_{st}$,  and the
direct image $ z_*(J): = dz \circ J \circ dz^{-1}$ satisfies $\vert\vert z_*(J) - J_{st}
\vert\vert_{C^m(\bar {\B})} \leq \lambda_0$.
\end{lemma}
\proof There exists a diffeomorphism $z$ from a neighborhood $U'$ of
$p \in M$ onto $\B$ satisfying $z(p) = 0$; after an additional linear change of coordinates one can achieve  $dz(p) \circ J(p)
\circ dz^{-1}(0) = J_{st}$ (this is a classical linear algebra). For $\lambda > 0$ consider the dilation
$d_{\lambda}: t \mapsto \lambda^{-1}t$ in $\R^{2n}$ and the composition
$z_{\lambda} = d_{\lambda} \circ z$. Then $\lim_{\lambda \rightarrow
0} \vert\vert (z_{\lambda})_{*}(J) - J_{st} \vert\vert_{C^m(\bar
{\B})} = 0$ for every  $m \geq 0$. Setting $U = z^{-1}_{\lambda}(\B)$ for
$\lambda > 0$ small enough, we obtain the desired statement. 
In what follows we often denote the structure $z^*(J)$ again by $J$ viewing it as a local representation of $J$ in the coordinate system $(z)$.

Recall that an almost complex structure $J$ is called {\it integrable} if $(M,J)$ is locally biholomorphic in a neighborhood of each point to an open subset of $(\C^n,J_{st})$. In the case of complex dimension $> 1$  integrable almost complex structures form a highly special subclass in the space of all almost complex structures on $M$.

\bigskip

In this paper we deal with  standard  classes of real submanifolds. A submanifold $E$ of an almost complex $n$-dimensional $(M,J)$ is called 
{\it totally real} if at every point $p \in E$ the tangent space $T_pE$  does not contain non-trivial complex vectors that is $T_pE \cap JT_pE = \{ 0 \}$. This is well-known that the (real)  dimension of a totally real submanifold of $M$ is not bigger than $n$; we will consider in this paper only $n$-dimensional totally real submanifolds that is the case of maximal dimension. A real submanifold $N$ of $(M,J)$ is called {\it generic} if the complex span of $T_pN$ is equal to the whole  $T_pM$ for each point $p \in N$. A real $n$-dimensional submanifold of $(M,J)$ is generic if and only if it is totally real.  
 
\begin{lemma}
\label{lemma2}
Let $N$ be a generic $(n+d)$-dimensional ($d > 0$) submanifold of an almost complex $n$-dimensional manifold $(M,J)$. Suppose that $K$ is a subset of $N$ of non-zero Hausdorff $(n+d)$-measure. Then 
there exists a (local) foliation of $N$ into a family $(E_s), s \in \R^{d}$ totally real $n$-dimensional submanifolds such that the intersection $K \cap E_s$ has a non-zero Hausdorff $n$-measure for 
each $s$ from some subset of non-zero Lebesgue measure in $\R^{d}$.
\end{lemma}
Here the Hausdorff measure is defined with respect to any Riemannian metric on $M$; the assumption that $K$ has a positive $n$-measure is independent on a choice of such metric.
\proof Let $p$ be a point of $M$ such that $K$ has a non-zero measure in each neighborhood of $p$. Choose local coordinates $z$ near  $p$ such that $p = 0$ and  $J(0) = J_{st}$.
After a $\C$-linear change of coordinates $N = \{ x_j + o(\vert z \vert) = 0, j = n-d+1,...,n\}$. After a local diffeomorphism with the identical linear part at $0$ we obtain that $N = \R^d(x_1,...,x_d) \times i\R^n(y)$. In the new coordinates the condition $J(0) = J_{s}$ still holds and every slice $E_s = \{ z \in N: x_1 = s_1,...,x_d= s_d\}$ is totally real. Now we conclude by the Fubini theorem.

\bigskip

A totally real manifold $E$ can be defined as 
 \begin{eqnarray}
 \label{edge}
 E = \{ p \in M: \rho_j(p) = 0 \}
 \end{eqnarray}
 where $\rho_j:M   \to  \R$  are smooth functions with non-vanishing gradients.  The condition of total reality means that for every $p \in E$ the $J$-complex linear parts of the differentials $d\rho_j$ are (complex) linearly independent.

 A subdomain

\begin{eqnarray}
\label{wedge}
 W = \{ p \in M: \rho_j < 0, j = 1,...,n \}.
 \end{eqnarray}t
 is called {\it the wedge with the edge} $E$.

 \subsection{Pseudoholomorphic discs} Let $(M,J)$ be an almost complex manifold of dimension $n > 1$. For a "generic" choice of an almost complex structure, any holomorphic (even locally) function on $M$ is constant. Similarly, $M$ does not admit non-trivial $J$-complex submanifolds (that is, with tangent spaces invariant with respect to  $J$) of complex dimension $> 1$. The only (but fundamentally important) exception arises in the case of pseudoholomorphic curves i.e. $J$-complex submanifolds of 
 complex dimension 1: they always exist locally.

Usually  pseudoholomorphic curves arise in connection with   solutions  $f$ of (\ref{CRglobal}) in the special case   where $M'$ has the complex dimension 1. These holomorphic maps are called $J$-complex (or $J$-holomorphic or {\it pseudoholomorphic} ) curves. Note that we view here the curves as maps i.e. we consider parametrized curves.
We use the notation  $\D = \{ \zeta \in \C: \vert \zeta \vert < 1 \}$ for  the
unit disc in $\C$ always assuming that it is equipped with the standard complex structure   $J_{\st}$. If in the equations (\ref{CRglobal})  we have $M' = \D$  we  call such a map $f$ a $J$-{\it complex  disc} or a
 {\it pseudoholomorphic disc} or just a  holomorphic disc
when the structure  $J$ is fixed. 

A fundamental fact is that  pseudoholomorphic discs always exist in a suitable neighborhood of any point of $M$; this is the classical Nijenhuis-Woolf theorem (see \cite{Aud}). Here it is convenient to rewrite the equations (\ref{CRglobal}) in local coordinates  similarly to the complex version of the usual Cauchy-Riemann equations.

 Everything will be local, so (as above) we are in a neighborhood $\Omega$ of $0$ in $\C^n$ with the standard complex coordinates $z = (z_1,...,z_n)$. We assume that $J$ is an almost complex structure defined on $\Omega$ and $J(0) = J_{st}$. Let 
$$z:\D \to \Omega,$$ 
$$z : \zeta \mapsto z(\zeta)$$ 
be a $J$-complex disc. Setting $\zeta = \xi + i\eta$ we write (\ref{CRglobal}) in the form $z_\eta = J(Z) Z_\xi$. This equation can be in turn written as

\begin{eqnarray}
\label{holomorphy}
z_{\bar\zeta} - A(z)\bar z_{\bar\zeta} = 0,\quad
\zeta\in\D.
\end{eqnarray}
Here a smooth map $A: \Omega \to Mat(n,\C)$ is defined by the equality $L(z) v = A \overline v$ for any vector $v \in \C^n$ and $L$ is an $\R$-linear map defined by $L = (J_{st} + J)^{-1}(J_{st} - J)$. It is easy to check that the condition $J^2 = -Id$ is equivalent to the fact that $L$ is $\overline\C$-linear. The matrix $A(z)$ is called {\it the complex matrix} of $J$ in the local coordinates $z$. Locally the correspondence between $A$ and $J$ is one-to-one. Note that the condition $J(0) = J_{st}$ means that $A(0) = 0$.

If $z'$ are other local coordinates and $A'$ is the corresponding complex matrix of $J'$, then, as it is easy to check, we have the following transformation rule:

\begin{eqnarray}
\label{CompMat}
A' = ({z'}_zA  + { z'}_{\overline z})({\overline z'}_{\overline z} + {\overline z'}_{ z}A)^{-1}
\end{eqnarray}
(see \cite{SuTu}).

Recall that for any complex function $f \in C^r(\D)$, $r >  0$,  {\it the Cauchy-Green transform} is defined by

\begin{eqnarray}
\label{CauchyGreen}
Tf(\zeta) = \frac{1}{2 \pi i} \int\int_{\D} \frac{f(\omega)d\omega \wedge d\overline\omega}{\omega - \zeta}
\end{eqnarray}
This is classical that the operator $T$ has the following properties:
\begin{itemize}
\item[(i)] $T: C^r(\D) \to C^{r+1}(\D)$ is a bounded linear operator for every non-integer $r > 0$ (similar property holds in the Sobolev scales). Here we use the usual H\"older norm on the space $C^r(\D)$.
\item[(ii)] $(Tf)_{\overline\zeta} = f$ i.e. $T$ solves the $\overline\partial$-equation in the unit disc. 
\item[(iii)] the function $Tf$ is holomorphic on $\C \setminus \overline\D$.
\end{itemize}
Fix a real non-integer $r > 1$. Let $z: \D \to \C^n$ be a $J$-complex disc. 
Since  the operator
$$\Psi_{J}: z \longrightarrow w =  z + TA(z) \overline {z}_{\overline \zeta} $$
takes the space   $C^{r}(\overline{\mathbb D})$  into itself,  we can write   the
equation (\ref{CRglobal}) in the form 
$[\Psi_J(z)]_{\overline \zeta}  = 0$. Thus, the disc $z$ is $J$-holomorphic if
and only if the map $\Psi_{J}(z):\mathbb D \longrightarrow \C^n$ is
$J_{st}$-holomorphic.
When the norm of $A$  is small enough (which is assured  by Lemma \ref{lemma1}),
then  by the implicit function theorem the operator    $\Psi_J$
is invertible  and we obtain a bijective
correspondence between  $J$-holomorphic discs and usual
holomorphic discs. This easily implies the existence of a $J$-holomorphic disc
in a given tangent direction through a given point of $M$, as well as  a smooth dependence of such a
disc  on a deformation of a point or a tangent vector, or on an almost complex structure; this also establishes  the interior elliptic regularity of discs. 

Note that global pseudoholomorphic discs (that is discs which are not contained in a small coordinate neighborhood) also have similar properties. The proof requires a considerably more subtil analysis of the integral operator $\Psi_J$, see \cite{SuTu2}.

\bigskip

Let $(M,J)$ be an almost complex manifold and $E \subset M$ be a real submanifold of $M$. 
Suppose that a $J$-complex disc $f:\D \to M$ is  continuous on $\overline\D$.  We some abuse of terminology, we also call the image $f(\D)$  simply by a disc and twe call he image $f(b\D)$ by  the boundary of a disc. If  $f(b\D) \subset E$, then we say that (the boundary of ) the disc  $f$ is {\it glued} or {\it attached} to $E$ or simply 
that $f$ is attached to $E$. Sometimes such maps are called {\it Bishop discs} for $E$ and we employ this terminology. Of course, if $p$ is a point of $E$, 
then the constant map $f \equiv p$ always satisfies this definition.

\section{Plurisubharmonic functions on almost complex manifolds} 
This section discusses some basic properties of plurisubharmonic functions on almost complex manifolds.

\subsection{Basic definitions}  Let  $u$ be a real $C^2$ function on an open subset $\Omega$ of an almost complex manifold  $(M,J)$. Denote by $J^*du$ the
 differential form acting on a vector field $X$ by $J^*du(X):= du(JX)$. Given point $p \in M$ and a tangent vector $V \in T_p(M)$ consider  a smooth vector field $X$ in a
neighborhood of $p$ satisfying $X(p) = V$. 
The value of the {\it complex Hessian} ( or  the  Levi form )   of $u$ with respect to $J$ at $p$ and $V$ is defined by $H(u)(p,V):= -(dJ^* du)_p(X,JX)$.  This definition is independent of
the choice of a vector field $X$. For instance, if $J = J_{st}$ in $\C$, then
$-dJ^*du = \Delta u d\xi \wedge d\eta$; here $\Delta$ denotes the Laplacian. In
particular, $H_{J_{st}}(u)(0,\frac{\partial}{\partial \xi}) = \Delta u(0)$.

Recall some  basic properties of the complex Hessian (see for instance, \cite{DiSu}):

\begin{lemma}
\label{pro1}
Consider  a real function $u$  of class $C^2$ in a neighborhood of a point $p \in M$.
\begin{itemize}
\item[(i)] Let $F: (M',J') \longrightarrow (M, J)$ be a $(J',
  J)$-holomorphic map, $F(p') = p$. For each vector $V' \in T_{p'}(M')$ we have
$H_{J'}(u \circ F)(p',V') = H_{J}(u)(p,dF(p)(V'))$.
\item[(ii)] If $f:\D \longrightarrow M$ is a $J$-complex disc satisfying
  $f(0) = p$, and $df(0)(\frac{\partial}{\partial \xi}) = V \in T_p(M)$ , then $H_J(u)(p,V) = \Delta (u \circ f) (0)$.
\end{itemize}
\end{lemma}
 Property (i) expresses the holomorphic invariance of the complex Hessian. Property (ii) is often useful in order to compute the complex Hessian  on a given 
tangent vector $V$.

 Let $\Omega$ be a domain $M$. An upper semicontinuous function $u: \Omega \to [-\infty,+\infty[$ on $(M,J)$ is
{\it $J$-plurisubharmonic} (psh) if for every $J$-complex disc $f:\D \to \Omega$ the composition $u \circ f$ is a subharmonic function on $\D$. Of course, this definition make sense because there are plenty of pseudoholomorphic discs  in a neighborhood of each point of an almost complex manifold.

By Proposition
\ref{pro1}, a $C^2$ function $u$ is psh on $\Omega$ if and only if it has    a 
positive semi-definite complex Hessian on $\Omega$ i.e. $H_J(u)(p,V) \geq 0$  for any $ p \in \Omega $ and $V \in T_p(M)$. 
A real $C^2$ function $u:\Omega \to \R$ is called {\it strictly $J$-plurisubharmonic} on $\Omega$, if $H_J(u)(p,V) > 0$ for each $p \in M$ and $V \in T_p(M) \backslash \{ 0\}$. Obviously, these notions   are local: an upper semicontinuous (resp. of class $C^2$) function on $\Omega$ is  $J$-plurisubharmonic (resp. strictly) on $\Omega$ if and only if it is $J$-plurisubharmonic (resp. strictly) in some open neighborhood of each point of $\Omega$. In what follows we often write "plurisubharmonic" instead of "$J$-plurisubharmonic" when an almost comple structure $J$ is prescribed.

A useful  observation is that the Levi form of a function $u$ at
a point $p$ in
an almost complex manifold $(M,J)$ coincides with the Levi form with respect
to the standard structure $J_{st}$ of $\R^{2n}$ if {\it suitable} local
coordinates near $p$ are choosen. Let us explain how to construct these adapted
coordinate systems.

As above, choosing local coordinates near $p$ we may identify a neighborhood
of $p$ with a neighborhood of the origin and assume that $J$-holomorphic discs
are solutions of (\ref{holomorphy}).

\begin{lemma}
\label{normalization}
There exists a  second order polynomial local diffeomorphism $\Phi$ fixing the
origin and with linear part equal to the identity such that in the new coordinates
 the complex matrix  $A$  of $J$ (that is $A$ from the equation (\ref{holomorphy})) satisfies
\begin{eqnarray}
\label{norm}
A(0) = 0, A_{z}(0) = 0
\end{eqnarray}
\end{lemma}
Thus, by a suitable local change of coordinates one can remove the 
terms linear in $z$ in the matrix $A$. We stress that in general it is impossible to
get rid of  first order terms containing  $\overline z$ since this would
impose a restriction on the Nijenhuis tensor $J$ at the origin.

I have learned this result  from  unpublished E.Chirka's notes; see .\cite{DiSu} for the proof. In
\cite{SuTu} it is shown that, in an almost complex manifold of (complex) dimension 2,
 a similar normalization is possible along a given embedded $J$-holomorphic disc.
 
 As an example consider a function $u(z) = \parallel z \parallel^2$ (we use the Euclidean norm) in the adapted coordinates (\ref{norm}). We conclude that this function is strictly $J$-plurisubharmonic 
  near the origin. In particular, each almost complex manifld admits plenty strictly $J$-psh functions locally. 
 
 As another typical consequence consider a totally real manifold $E$ defined by (\ref{edge}). Then the function $u = \sum_{j=1}^n \rho^2$ is strictly $J$-plurisubharmonic  in a neighborhood of $E$. Indeed, it suffices to choose local coordinates near $p \in M$ according to Lemma \ref{normalization}. This reduces the verification to the well-known case of $J_{st}$.

 \subsection{Envelopes of plurisubharmonic functions} It follows  from definitions that  many of the elementary properties of plurisubharmonic functions can be directly transferred to the almost complex case. We mention here, for example, the maximum principle as well as the fundamental fact that the plurisubharmonicity is a local property: a function is plurisubharmonic on $M$ if and only if it is plurisubharmonic in an open  neighbourhood of every point of $M$.

Let $(M,J)$ be an almost complex manifold of complex dimension $n$. Fix a Riemannian metric on $M$;
all norms and distances will be considered with respect to this metric. Of course, the results are independent on its choice.

 As another typical example we recall here a construction of an envelope of a family of plurisubharmonic functions following S.Bu - W.Schaschermayer \cite{BuScha} (in the almost complex case this construction was used in \cite{CoGaSu}).

Let 
\begin{eqnarray}
\label{average}
P_0\phi = \frac{1}{2\pi i} \int_{b\D} \phi(\omega) \frac{d\omega}{\omega}
\end{eqnarray}
denotes the average of a real function $\phi$ over $b\D$,

\begin{lemma}
\label{lemma3}
Let $v$ be an upper semicontinuous function on an almost complex manifold 
$(M,J)$. Consider the sequence $(v_n)$ defined as follows: $v_0 = v$ and
for $m \geq1$, for $p \in M$,
$$
v_m(p) = \inf_f \,\, P_0 (v_{n-1} \circ f),
$$
where $\inf$ is taken over all $J$-complex discs $f : \D \to M$ such that $f(0) = p$, $f$ is of class 
$C^r(\overline\D)$ with some (fixed) non-integer $r > 1$ and $f(\overline\D) \subset M$. Then the sequence $(v_n)$
decreases pointwise to the largest $J$-plurisubharmonic function $DE[v]$ on $M$
bounded from above by $v$.
\end{lemma}

\proof We proceed in several steps. Clearly, every $v_m$ is correctly defined. 

{\it Step 1.} The sequence $(v_m)$ decreases.
Indeed, for every $p$, the constant disc $f^0(\zeta) \equiv p$ is
$J$-complex so
$$v_m(p) = \inf_f \,\, P_0(v_{m-1} \circ f) \leq
P_0 (v_{m-1} (p))  = v_{m-1}(p).
$$

{\it Step 2.}   The function $DE[v]$ is upper semicontinuous. 
We proceed the proof by induction on $m$. For $m = 0$ the statement is
correct. Suppose that the function $v_{m-1}$ is upper semicontinuous. Let $(p_k) \subset M$ be a
sequence of points converging to $p_0 \in M$.

It follows from  \cite{SuTu2} that the following holds. Let $f:\D \to M$ be a $J$-complex disc
  of class  $C^r(\D)$ such that $f(0) = p_0$ and $f(\overline \D) \subset M$. Then there
exists a sequence of $J$-complex discs $f_k : \D \to M$, of class $C^r(\D)$ such that $f_k(\overline \D) \subset M$, $f_k(0) = p_k$
for every $k$ and $f_k \longrightarrow f$ in  $C^r(\D)$.

\vskip 0,1cm
Consider a compact set $K$ containing   the union $f(\overline\D) \cup (\cup_k
f_k(\overline \D))$. Since $v_{m-1}$ is an upper semicontinuous function, it is upper bounded on $K$  by a
constant $C$, and
$$
(v_{m-1}\circ f)(\zeta) \geq \lim \sup_{k \to \infty}
(v_{m-1} \circ f_k)(\zeta),\,\,\, \zeta \in \D.
$$

So by the Fatou lemma
$$
P_0 (v_{m-1} \circ f) \geq \lim \sup_{k
\to \infty} P_0(v_{m-1} \circ f_k) \geq \lim
\sup_{k \longrightarrow \infty} v_m(p_k).
$$
This implies that
$$
v_m(p_0) = \inf_f \,\, P_0(v_{m-1} \circ f)\geq \lim \sup_{k
\to \infty} v_m(p_k)
$$
which proves  the upper semicontinuity of $v_m$.

Therefore, the function $DE[v]$ also is upper semicontinuous as a decreasing
limit of upper semicontinuous functions.

{\it Step 3.} We prove by induction that for any $J$-plurisubharmonic function
$\phi$ satisfying $\phi
\leq v$ we have $ \phi \leq v_n$ for any $n$.  This is true for $m = 0$.
Suppose that $\phi \leq v_{m-1}$. Fix an arbitrary point $p_0 \in M$.
For every $J$-complex disc $f \in C^r(\D)$ and satisfying $f(0) =p_0$,
$f(\overline \D) \subset M$, we have~:
$$
\phi(z_0) \leq
P_0(\phi \circ f ) \leq P_0 (v_{m-1} \circ f) .
$$
 Hence
$v_m(p_0) \ge \phi(p_0)$.

{\it Step 4.} We show that the restriction of $DE[v]$ on a
$J$-complex disc is subharmonic. Given $p_0$ and $f$ as above, we have
$$
DE[v](p_0) = \lim_{m \to \infty} v_m(p_0) \leq \lim_{
\to \infty} P_0 (v_{m-1} \circ f)  = P_0 (v \circ f)
$$
by the Beppo Levi theorem. This concludes the proof of lemma. 

We call the function $DE[v]$  {\it the disc envelope} of $v$. As a simple application, consider 
any family $(u_\alpha)$ of plurisubharmonic functions on $(M,J)$ and the function 
$u= \sup_\alpha \, u_\alpha$. In general, $u$ does not need to be upper semicontinuous, so we consider its {\it upper regularization} 
$$u^*(p) = \lim_{\varepsilon \to 0+} \inf_{dist(q,p) \le \varepsilon} u(q)$$
It is classical that $u^*$ is the smallest upper semicontinuous function satisfying $u \le u^*$.
In order to prove that $u^*$ is plurisubharmonic on $M$, consider the disc envelope $DE[u^*]$. It follows from Lemma \ref{lemma3} that $u_\alpha \le DE[u^*]$ for all $\alpha$ that is $u \le DE[u^*] \le u^*$. Hence $DE[u^*] = u^*$ and $u^*$ is plurisubharmonic. Usually $u^*$ is called the {\it $\sup$-envelope} of the family $(u_\alpha)$. 

Note that the usual (for $M = \C^n$) proofs of plurisubharmonicity of $u^*$ do not go through directly in the almost complex case since they are based on regularization of plurisubharmonic functions by convolution (see for example \cite{De}). This argument is not available in the general almost complex case because the Cauchy-Riemann equations (\ref{holomorphy}) are only quasi-linear and not linear.

We point out here a difference of the above construction of the disc envelope and the argument of \cite{BuScha}. In \cite{BuScha} only linear complex discs are used. Of course, this does not make sense in the almost complex case and we need to consider all pseudoholomorphic discs. As a consequence, the set of discs under consideration is much larger and from this point of view, we are closer to the construction of the disc envelope introduced by E.Poletsky \cite{Po}. He proved (in the case of $\C^n$) that the iteration process used in the proof of Lemma \ref{lemma3} stopes already on the first step, that is 
$v_1 =v_2=....= DE[v]$. His result was extended to the case of complex manifolds by 
F.Larusson - R.Sigurdsson \cite{La,Ro2} and J.-P.Rosay \cite{Ro2}.  To the best of my knowledge, it remains an open question if this is also true for  almost complex manifolds in any dimension;  the 
case of dimension 2 is settled by U.Kuzman \cite{Ku}.

\subsection{Plurisuperharmonic measure}  
A function $v$ is called {\it plurisuperharmonic} on $M$ if the function $-v$ is plurisubharmonic on $M$.

Let $\Omega$ be a smoothly bounded domain in $M$ with the boundary $b\Omega$. 
For $\alpha \ge 1$ and $q \in b\Omega$ set 
$A_\alpha(q) = \{ p \in \Omega: dist(p,q) \le \alpha d_q(p) \}$. Here $dist$ denotes the distance on $M$ and $d_q(p)$ denotes the distance from $p$ to the tangent plane $T_q(b\Omega)$ to $b\Omega$ at $q$. In the case where $M = \C^n$ with the standard Euclidean distance, this is the intersection of $\Omega$ with a cone with vertex at $w$. In the general case $A_\alpha(q)$ is a region of $\Omega$ which approaches  $b\Omega$ non-tangentially at $q$.

Let $u$ be a plurisuperharmonic function on $\Omega$. Denote by $u_*$ its {\it non-tangential lower boundary extension}  which is defined as 
$$u_*(q) = u(q),\,\,\, q \in \Omega,$$  
$$u_*(q) =  \inf_{\alpha > 1} (\lim \inf_{A_\alpha(q) \ni p \to q} u(p)), \,\,\, q \in b\Omega.$$

Let $K$ be a compact subset of $\overline\Omega$. Denote by $P(K)$ the class of all functions $u$ plurisuperharmonic on $\Omega$ and such that $u(q) \ge 0$ for each $q \in \Omega$ and $u_*(q) \ge 1$ for each $q \in K$. {\it The plurisuperharmonic measure of $K$ with respect to $\Omega$} or simply {\it the $P$-measure} is the function
\begin{eqnarray}
\label{Pmeasure1}
\omega_*(p,K,\Omega) = \lim\inf_{q \to p} \omega(q,K,\Omega)
\end{eqnarray}
where
\begin{eqnarray}
\label{Pmeasure2}
\omega(q,K,\Omega) = \inf_{u \in P(K)} u(q)
\end{eqnarray}
Of course, in the one-variable case $\omega_*$ coincides with the usual harmonic measure.

\bigskip

Following \cite{Sa1}, consider some basic properties of the $P$-measure. 
Recall that a subset $E \subset \Omega$ is called {\it pluripolar} if there exists a plurisuperharmonic on $\Omega$ function $u$ non identically equal to $+\infty$ and such that $u\vert E = + \infty$. It follows from the results of F. R.Harvey - H.B.Lawson \cite{HL} that a pluripolar set is of measure zero.

\begin{prop}
\label{PMProp1}
We have the following
\begin{itemize}
\item[(i)] $0 \le \omega_*(p,K,\Omega) \le 1$ for every $p \in \Omega$.
\item[(ii)] The function $p \mapsto \omega_*(p,K,\Omega)$ is plurisuperharmonic on $\Omega$.
\item[(iii)] If $\Omega_1 \subset \Omega_2$ and $K_1 \subset K_2$, then $\omega_*(p,K_1,\Omega_1) \le \omega_*(p,K_2,\Omega_2)$.
\item[(iv)] If $\omega_*(p^0,K,\Omega) = 0$ for some $p^0 \in \Omega$, then $\omega_*(p,K,\Omega) = 0$ for all $p \in \Omega$.
\item[(v)] Let $K \subset \overline\Omega$ is a pluripolar subset of some open neighborhood $\tilde\Omega$ of $\overline\Omega$. Then $\omega_*(p,K,\Omega) = 0$ for all $p \in \Omega$.
\end{itemize}
\end{prop}
The property (i) is obvious (the second inequality follows because the constant function $u = 1$ belongs to $P(K)$). The property (ii) is classical. The property (iii) is obvious; (iv) follows by the maximum principle. For (v), suppose that there exists a function $u$ plurisuperharmonic on $\tilde\Omega$ and which is not equal to $+\infty$ identically such that $u \vert K = + \infty$. One can assume that $u \ge 0$ on $\Omega$. For $m=1,2,...$ the function 
$$v_m(p) = \min (u(p)/m, 1)$$
is superharmonic on $\Omega$ and belongs to the class $P(K)$. Hence 
$$\omega(p,K,\Omega) \le v(p) = \lim_{m \to \infty} v_m(p)$$
But $v(p) = 0$ when $u(p) \neq +\infty$ and $v(p) =  1$ when $u(p) = +\infty$. Therefore the lower regularization $\omega_*(p,K,\Omega)$ vanishes identically.

\bigskip

As the first application let us prove the following version of two constants theorem. Denote by $u^*$ the {\it upper non-tangent boundary extension}  of a plurisubharmonic function $u$ that is 

$$u^*(q) = u(q),\,\,\, q \in \Omega,$$  
$$u^*(q) =  \sup_{\alpha > 1} (\lim \sup_{A_\alpha(q) \ni p \to q} u(p)), \,\,\, q \in b\Omega.$$

Similarly to \cite{Sa1}, we have

\begin{prop}
\label{2Const1}
Let $u$ be a plurisubharmonic function bounded above by $C$ on a smoothly bounded domain $\Omega$ in $(M,J)$.  Suppose that for some compact subset $K \subset \overline\Omega$ we have
$u^*(p) \le c < C$ for every $p \in K$. Then
$$u(p) \le c\omega_*(p,K,\Omega) + C(1 - \omega_*(p,K,\Omega))$$
\end{prop}
For the proof it suffices to note that $\omega_*(p,K,\Omega) \le (C-u(p))/(C-c)$ 
because the function in the right hand is of class  $P(K)$.

\section{Construction of complex discs}

This section presents our main technical tool. We  fill a wedge $W$ with a totally real edge $E$ by a family of complex  discs  glued to the edge $E$ along the upper semi-circle. We apply the approach developed in \cite{Su1} which requires some refinement suitable for our goals.

We will proceed in several steps.

\bigskip

(a)   First consider the model case where $M = \C^n$ with $J = J_{st}$ and $E = i\R^n = \{ x_j = 0, j = 1,...,n \}$. Denote by $W$ the standard wedge $W = \{ z= x +iy: x_j < 0, j=1,....,n \}$.

Consider the family of linear complex maps

\begin{eqnarray}
\label{disc1}
l: (c,t,\zeta) \mapsto (\zeta, \zeta t + ic )
\end{eqnarray}
Here $\zeta \in \C$; the variables  $c = (c_2,...,c_{n}) \in \R^{n-1}$ and $t\in \R^{n-1}_+= \{t = (t_2,...,t_n) \in \R^{n-1}: t_j > 0\}$) are viewed as parameters.
Denote by $V$ the wedge $V = \R^{n-1} \times \R^{n-1}_+$. Also let $\Pi = \{ \Re \zeta < 0 \}$ be the left half-plane; its boundary $b\Pi$ coincides with the imaginary axis $i\R$. The following properties of the above  family are easy to check:
\begin{itemize}
\item[(a1)] the images  $l(c,t)(b\Pi)$  form a family of  real lines in  $i\R^n = E$ . For every fixed 
$t \in \R^{n-1}_+$  these lines are disjoint and  
$$\cup_{c \in \R^{n-1}} l(c,t)(b\Pi) = E.$$
In other words, for every $t$ this family (depending on the parameter $c$) forms a foliation of $E$ by 
parallel lines. 
\item[(a2)] one has
$$\cup_{(c,t) \in V} l(c,t)(\Pi) = W.$$
\item[(a3)] For every fixed $t \in \R^{n-1}_+$, one has
$$\cup_{c \in \R^{n-1}} l(c,t)(\Pi) = E_t = \{ z \in \C^n: \Re(z_j - t_j z_1) = 0, j=2,...,n \} \cap W$$
and the union is disjoint.
Every $E_t$ is a real linear $(n+1)$-dimensional half-space contained in $W$ and $bE_t = E$.
\item[(a4)] the family $(E_t)$, $t \in \R^{n-1}_+$ is  disjoint in $W$ and its union coincides with $W$.
\end{itemize} 

Let $K \subset E$ be a compact subset of non-zero Hausdorff $n$-measure (one can consider it with respect to the standard metric). Consider the set $\Sigma_t$ of $c \in \R^{n-1}$ such that the real line $l(c,t)(b\Pi)$ intersects $K$ in a subset of non-zero $1$-measure. 
It follows by (a1) and the Fubini theorem that for every $t \in \R^{n-1}_+$ ,  the set $\Sigma_t$ has a non-zero $(n-1)$-measure. Again by the Fubini theorem and (a3), the set 
$\cup_{c \in \Sigma_t} l(c,t)(\Pi)$ is a subset of $E_t$ of non-zero $(n+1)$ - measure. Finally, by (a4) we obtain 

\begin{itemize}
\item[(a5)] $\cup_{t \in \R^{n-1}_+} \cup_{c \in \Sigma_t} l(c,t)(\Pi)$ is a subset of $W$ of non-zero $2n$-measure.
\end{itemize}
In what follows we will use these properties locally that is in a neighborhood of the origin. It is convenient to reparametrize the family of complex half-lines $l(c,t)$ by complex discs.

\bigskip

We represent the family of discs (\ref{disc1}) (after suitable reparametrization) as a general solution of an integral equation.

Let 
\begin{eqnarray}
\label{Schwarz}
S \phi(\zeta) = \frac{1}{2\pi i} \int_{b\D} \frac{\omega + \zeta}{\omega - \zeta }\phi(\omega) \frac{d\omega}{\omega}
\end{eqnarray}
denotes the Schwarz integral. In terms of the Cauchy transform 

\begin{eqnarray}
\label{Cauchy}
Kf(\zeta) = \frac{1}{2 \pi i} \int_{b\D} \frac{f(\omega)d\omega}{\omega - \zeta}
\end{eqnarray}
we have the following relation: $S = 2K - P_0$. As a consequence the boundary properties of the Schwarz integral are the same as the classical properties of the Cauchy integral.

For a {\it non-integer} $r > 1$ consider the Banach spaces $C^r(b\D)$ and $C^r(\D)$ (with the usual H\"older norm). This is classical that $K$ and $S$ are bounded linear mappings in these classes of functions. For a real function $\phi \in C^r(b\D)$ the Schwarz integral $S\phi$ is a function of class $C^r(\D)$ holomorphic in $\D$; the trace of its real part on the boundary coincides with $\phi$ and its imaginary part vanishes at the origin. In particular,  every holomorphic function $f \in C^r(\D)$ satisfies the Schwarz formula $f = S\Re f + iP_0f$.

We are going to fill $W$  by complex discs glued to $i\R^n$ along the (closed) upper 
semi-circle $b\D^+ = \{ e^{i\theta}: \theta \in [0, \pi] \}$; let also $b\D^-:= b\D \setminus b\D^+$.

Fix a smooth real function $\phi: b\D \to \R$ such that $\phi \vert b\D^+ = 0$ and $\phi \vert b\D^- <  0$.

Consider now  a real $2n$-parameter  family of holomorphic discs $z^0 = (z_1^0,...,z_n^0) : \D \to \C^n$ with components

\begin{eqnarray}
\label{BP1}
z_j^0(c,t)(\zeta) = x_j(\zeta) + iy_j(\zeta) = t_j S\phi(\zeta) + ic_j, j= 1,...,n
\end{eqnarray}
Here $t_j > 0$ and $c_j \in \R$ are  parameters, $t= (t_2,...,t_n)$, $c= (c_2,...,c_n)$; in (\ref{BP1})
we formally set $t_1 = 1$ and $c_1= 0$.

Obviously, every $z^0(c,t)(\D)$ is a subset of $l(c,t)(\Pi)$ and $z^0(b\D^+) = l(c,t)(b\Pi)$. Thus, the family $z^0(c,t)$ is a (local) biholomorphic reparametrization of the family $l(c,t)$. As a consequence, the properties (a1)-(a5) also hold for the family $z^0(c,t)$.
Notice also the following obvious properties of this family:

\begin{itemize}
\item[(a6)] for every $j$ one has $x_j \vert b\D^+ = 0$ and $x_j(\zeta) < 0$ when $\zeta \in \D$ (by the maximum principle for harmonic functions). 
\item[(a7)] the evaluation map $Ev_0: (c,t,\zeta) \mapsto z^0(c,t)(\zeta)$ is one-to-one from $V \times \D$ to $W$.
\end{itemize}

Now we construct an  analog of this family in the general case.

\bigskip

(b)  In order to write an integral equation defining a required family of discs we need to employ an analog of the Schwarz formula and to choose suitable local coordinates.

 We have the following Green-Schwarz formula (see the proof, for example, in  \cite{SuTu}, although, of course, it can be found in the vaste list of classical works). Let $f = \phi + i\psi : \overline\D \to \C$  be a function of class $C^{r}(\D)$. Then for each $\zeta \in \overline\D$ one has
\begin{eqnarray}
\label{GreenSchwarz}
f(\zeta) = S \phi(\zeta) + i P_0 \psi + Tf_{\overline\zeta}(\zeta) - \overline{Tf_{\overline\zeta}(1/\overline\zeta)}
\end{eqnarray}
Here we use the integral operators (\ref{CauchyGreen}), (\ref{average}) and  (\ref{Schwarz}). Note that because of the "symmetrization" the real part of the sum of two terms containing the Cauchy-Green operator $T$ vanishes on the unit circle. Notice also that the last term is holomorphic on $\D$.

 Now let $(M,J)$ be an almost complex manifold of complex dimension $n$ and $E$ be a totally real $n$-dimensional submanifold of $M$.  We assume that $E$ and $W$ are given by (\ref{edge}) and (\ref{wedge}) respectively.

First, according to Section 2 we choose local coordinates $z$ such that $p = 0$ and the complex matrix $A$ of $J$ satisfies (\ref{norm}). For every $\tau > 0$ small enough and $C > 0$ big enough the functions 
$$\tilde\rho_j:=  \rho_j - \tau \sum_{ k \neq j} \rho_k  + C\sum_{k=1}^n \rho_k^2$$
are strictly $J$-plurisubharmonic  in a neighborhood of the origin and the "truncated" wedge $W_\tau = \{ \tilde\rho_j < 0, j=1,...,n \}$ is contained in $W$. After a $\C$-linear (with respect to $J_{st}$) change of coordinates 
one can assume that $\tilde\rho_j = x_j + o(\vert z \vert)$. Consider now a local diffeomorphism 
$$\Phi: z = x_j + iy_j \mapsto z' = x_j' + iy_j' = \tilde\rho_j + iy_j$$
Then $\Phi(0)= 0$, $d\Phi(0) = 0$ and in the new coordinates $\tilde\rho_j = x_j$ (we drop the primes), $E = i\R^n$ and $W_\tau = \{ x_j < 0 , j= 1,...,n\}$. We keep the notation $J$ for the direct image 
$\Phi_*(J)$. Then in the choosen coordinates the complex matrix of $J$ still satisfies  $A(0) = 0$. Note also that the coordinate functions $x_j$ are strictly plurisubharmonic  for $J$. 

Finally, similarly to the proof of Lemma \ref{lemma1}, for $\lambda > 0$  consider the isotropic dilations $d_\lambda: z \mapsto \lambda^{-1}z$ and the direct images $J_\lambda:= (d_\lambda)_*(J)$. Denote by $A(z,\lambda)$ the complex matrix of $J_\lambda$.

 For $\lambda > 0$ small enough we are looking for the solutions $z: \D \to \C^n$ of Bishop's type integral equation

\begin{eqnarray}
\label{MainBishop}
z(\zeta) = h(z(\zeta), c,t,\lambda) 
\end{eqnarray}
with
$$h(z(\zeta), c,t,\lambda) = tS \phi(\zeta) + i c + TA(z,\lambda)\overline{ z}_{\overline\zeta}(\zeta) - \overline{T A(z,\lambda) \overline{z}_{\overline\zeta}(1/\overline\zeta)}$$
where $t = (t_2,...,t_n)$, $t_j > 0$ and $c \in \R^{n-1}$ as well as $\lambda$ are viewed as real parameters. Of course, we assume implicitely that all parameters are close to the origin.

Note that the first and the last terms in the right hand are holomorphic on $\D$. Therefore, any solution of (\ref{MainBishop}) satisfies the Cauchy-Riemann equations (\ref{holomorphy}) i.e. is a $J$-complex disc. Furthemore, $x_j(\zeta)$ vanishes on $b\D^+$ ( that is $z(b\D^+) \subset E$) and is negative on $b\D^-$. Since the function $z \mapsto x_j$ is strictly $J$-plurisubharmonic, by the maximum principle the image $z(\overline\D)$ is contained in $\overline W_\tau$.

The existence of solutions follows by the implicit function theorem. Note that for $\lambda = 0$ the equation ( \ref{MainBishop}) admits the solution ( \ref{BP1}). Consider the smooth map
of Banach spaces

$$H: C^r(\D) \times \R^{n-1} \times \R^{n-1} \times \R  \longrightarrow  C^r( \D)$$
$$ H: (z, c,t, \lambda)  \mapsto h(z(\zeta),c,t,\lambda)$$ 

Obviously the partial derivative of $H$ in $z$ vanishes: $(D_zH)(z^0,c,t,0) = 0$, where $z^0$ is a disc given by ( \ref{BP1}). By the implicit function theorem, for every $(c,t, \lambda)$ close enough to the origin the equation ( \ref{MainBishop}) admits a unique solution 
\begin{eqnarray}
\label{disc2}
(c,t,\zeta) \mapsto z(c,t)( \zeta)
\end{eqnarray}
 of class $C^r( \D)$, smoothly depending on parameters $(c,t)$  (as well as  $\lambda$, of course). 

Fixing $ \lambda > 0$, consider the smooth evaluation map
$$Ev_ \lambda: (c,t,\zeta)  \mapsto z(c,t)(\zeta)$$
which associates to each parameter a point of the corresponding disc. Note that for $ \lambda = 0$ we  obtain the linear mapping $Ev_0(c,t)(\zeta)$ appeared already  in (a7) for the model case of the standard structure. Indeed, in this case the family  (\ref{disc2}) coincides with the family (\ref{BP1})  (by the uniqueness of solutions assured by the implicit function theorem). Notice also that  when $t = 0$, for every $\lambda > 0$ and  every $c \in \R^{n-1}$     the equation (\ref{MainBishop}) has the unique solution $z^0(c,0)(\zeta) =(\zeta, ic)$ (cf. with (\ref{BP1})). Hence, $Ev_\lambda(c,0)(\zeta) = Ev_0(c,0)(\zeta)$.

By (a7) $Ev_0(V \times \D)$ coincides with $W_\tau$ in a neighborhood of the origin.  Furthermore, 
$Ev_\lambda(\{ (c,t): c \in \R^{n-1}, t = 0 \} \times b\D^+) = E$ for $\lambda > 0$. Also, for $\alpha > 0$ the "truncated" wedge 
$W_\alpha = \{ x_j - \alpha\sum_{k \neq j} x_k < 0 \}$ with the edge $E$ is contained in $W_\tau$. The faces of the boundary of $W_\alpha$ are transversal to the face of $W_\tau$. Since this property is stable under small perturbations,  we conclude that $W_\alpha \subset Ev_\lambda(V)$ for all $\lambda$ small enough. In terms of the initial defining functions $\rho_j$ we have $\{ z: \rho_j - \delta\sum_{k \neq j} \rho_k < 0 \} \subset W_\varepsilon$ when $\tau + \varepsilon < \delta$.

 Concerning the regularity of manifolds and almost complex structures, it suffices to require the class $C^r$ with real $r > 2$ and the argument goes through. We skip the details.

\bigskip

Fix $\delta > 0$. Since the properties of linear discs (a1)-(a5) are stable under small perturbations, the obtained family of discs admits the similar propeties. For reader's convenience we list them.
\begin{itemize}
\item[(b1)] the images  $z(c,t)(bD^+)$  form a family of  real curves in  $E$ . For every fixed 
$t \in \R^{n-1}_+$  these curves are disjoint and  
$$\cup_{c \in \R^{n-1}} z(c,t)(b\D^+) = E.$$
In other words, for every $t$ this family (depending on the parameter $t$) forms a foliation of $E$. 
Furthermore, every disc is contained in $W$.
\item[(b2)] one has
$$\cup_{(c,t) \in V} z(c,t)(\D) = W_\delta = \{ z: \rho_j - \delta\sum_{k \neq j} \rho_k < 0 \}.$$
\item[(b3)] For every fixed $t \in \R^{n-1}_+$, the union
$$E_t:= \cup_{c \in \R^{n-1}} z(c,t)(\D) \subset W_\delta$$ is a real  $(n+1)$-dimensional manifold with boundary  $bE_t = E$.
\item[(b4)] the family $(E_t)$, $t \in \R^{n-1}_+$ is  disjoint and its union coincides with $W_\delta$.
\end{itemize} 

 Similarly  (a5) we have

\begin{itemize}
\item[(b5)] Let $K \subset E$ be a compact subset of non-zero Hausdorff $n$-measure. The discs $z(c,t)$ whose boundaries intersect $K$ in a set of positive $1$-measure, fill a subset of $W_\delta$  of non-zero Hausdorff $2n$-measure.
\end{itemize}

Consider now the important special case where a totally real $n$-dimensional manifold  $E$ of the form (\ref{edge}) is contained in the boundary $b\Omega$ of a domain $\Omega$. We assume that $b\Omega$ is a smooth real hypersurface  defined in a neighborhood of a point $p \in E \subset b\Omega$ by $b\Omega = \{ \rho = 0 \}$ where $\rho$ is a smooth real function with non-vanishing gradient; one can assume also that $\Omega = \{ \rho < 0 \}$. Then we can always shoose a wedge $W$ of the form (\ref{wedge}) such that $W \subset \Omega$ and its faces $\{\rho_j = 0\}$ are 
transversal to $b\Omega$. Then each complex disc from the constructed above family also is 
transversal to $b\Omega$. More precisely, we have the following property

\begin{itemize}
\item[(b6)] Every disc from the  family (\ref{disc2}) is non-tangent to $b\Omega$ at every boundary point.
\end{itemize}

\section{Boundary behavior of plurisubharmonic functions }

Now we are able to prove our main results. 

\begin{thm}
\label{theo2}
Let $\Omega$ be a smoothly bounded domain in an almost complex n-dimensional manifold $(M,J)$.
Suppose that $E \subset b\Omega$ be a generic submanifold and $K \subset E$ be a  subset  with non-empty interior (with respect to $E$). Then $\omega_*(p,K,\Omega)$ does not vanish identically.
\end{thm}
Here the Hausdorff measure is considered with respect to any Riemannian metric on $M$; as we already have pointed out, the condition of non-vanishing of the Hausdorff measure of $K$ is independent of a choice of the Riemannian metric.

In view of Lemma \ref{lemma1}, without loss of generality one can assume that $E$ is totally real. 
Consider the family (\ref{disc2}) of discs for $E$ constructed in Section 4; recall that these discs are not tangent to the boundary $b\Omega$ in view of (b6). Since the construction is local and $K$ has a non-empty interior in $E$, one can assume $K= E$. Using the properties (b4) and (b6) we see that the discs fill an open subset $X$ of $\Omega$. Each disc $f$ intersect $K$ along $b\D^+$; hence for every  $\zeta \in \D$  with $\Im \zeta > 0$ one has
$$\omega(f(\zeta),K,\Omega) \ge \omega(\zeta,b\D^+,\D)  \ge c > 0$$
where $c$ is a  universal constant. We obtain that 
$\omega(p,K,\Omega) \ge c > 0$ for each $p \in X$. We conclude that 
$\omega_*(p,K,\Omega) > 0$ on $X$  that is does not vanish identically.

\bigskip

The next result is the following uniqueness principle.

\begin{thm}
\label{theo1}
Let $\Omega$ be a smoothly bounded domain in an almost complex n-dimensional manifold $(M,J)$.
Suppose that $E \subset b\Omega$ be a generic submanifold and $K \subset E$ be a compact subset of non-zero n-Hausdorff measure.  Assume that $u$ is an upper bounded  plurisubharmonic function in  $\Omega$ such that $u^*(p) = -\infty$ for each $p \in K$. Then $u \equiv - \infty$.
\end{thm}

Indeed, for every disc $f$ from the family (\ref{disc2})   the composition  $u \circ f$ is an upper bounded 
subharmonic function  in $\D$. Notice again that the boundary of every disc is transversal to $b\Omega$ (see (b6)).  Consider now the discs $f$ from these family intersecting $K$ along a subset $l_f \subset b\D$ of positive measure. For each disc $f$ and $\zeta \in \D$ one has
$$\omega(f(\zeta),K,\Omega) \ge \omega(\zeta,l_f,\D) > 0$$
The above discs fill a subset $X$ of $\Omega$ of positive $2n$-measure and with non-empty interior in $\Omega$ according to (b5). Applying to every $u \circ f$
the two constant theorem for subharmonic functions, we conclude that $u\circ f \equiv -\infty$. Since these discs fill a set of positive measure in $\Omega$, we conclude that $u \equiv = -\infty$. 

\bigskip

In particular, we obtain the following far reashing generalization of one of the results of J.-P.Rosay \cite{Ro1}.

\begin{cor}
Let $E$ be a totally real $n$-dimensional submanifold of an almost complex $n$-dimensional manifold $(M,J)$. Suppose 
that $K$ is a  closed subset of $E$ of non-zero Hausdorff $n$-measure. Then $E$ is not contained in a pluripolar set.
\end{cor}

 We need some additional properties of discs (\ref{disc2}) constructed in Section 4. We use the notation from that section.

\bigskip

(a) Once again we begin with the standard case  where $M = \C^n$ with $J = J_{st}$ and $E = i\R^n = \{ x_j = 0, j = 1,...,n \}$. Consider the family $l(c,t)$ of the form (\ref{BP1}) (that is  (\ref{disc1}) after a reparametrization); these complex maps attach the left half-plane $\Pi$ to $E$ along the imaginary axis. Our goal is to study more carefully those maps from this family whose boundary does not touch the origin in $i\R^n$. This means that $c_j \neq 0$ for at least one $j \in \{2,....,n\}$. Each  point $z = x + iy \in W = \{ x_j < 0, j=1,...,n \}$ belongs to the disc $l(c,t)$ with $t_j = x_1/x_j > 0$ and $c_j = y_j - y_1x_1/x_j$. For each $j=2,...,n$ consider a real smooth hypersurface $\Gamma_j = \{ z \in W: y_j x_j - y_1 x_1 = 0 \}$ in $W$. We obtain the following

\begin{itemize}
\item[(a8)] the discs (\ref{disc1}) whose boundary does not touch the origin fill an open dense subset  $W \setminus \cup_{j=2}^n \Gamma_j$
of $W$.
\end{itemize}
Now we pass to the general case.

\bigskip

(b) Assume that $E$ is contained in the boundary of a smoothly bounded domain $\Omega$ and  consider discs (\ref{disc2}) such that (b6) holds.
Since the property (a8) is stable under perturbation, its analog holds for the family of pseudoholomorphic discs (\ref{disc2}). 
We consider only non-constant  discs with $t \in \R^n_+$.
\begin{itemize}
\item[(b7)] the discs $z(c,t)$ of the family (\ref{disc2}) whose boundary does not touch a point $p \in E$ fill an open dense subset  $X$ of $W_\delta$. 
\item[(b8)] In particular, under the assumptions of (b6), the closure $\overline X$ of $X$ contains any 
non-tangential region $A_\alpha(p)$ in a domain $\Omega$.
\end{itemize}

The following property that  we need  is obvious as well:

\begin{itemize}
\item[(b9)] under the  assumptions of (b7), the arcs $z(c,t)(b\D^-)$ fill a compact subset $Y\in W_\delta$.
\end{itemize}

\bigskip

Our next result is inspired by the work of Y.Khurumov \cite{Khu}. Another motivation arises from the work by N.Levenberg - G.Martin - E.Poletsky \cite{LMP}. They proved that the direct analog of the Lindel\"ef-Chirka principle \cite{Ch} fails for upper bounded plurisubharmonic function in complex dimension $> 1$. The following results can be viewed as a partial analog of two constant theorem and Lilndel\"ef - Chirka principle for plurisubharmonic functions.

\begin{thm}
\label{Lindel1}
Let $\Omega$ be a smoothly bounded domain in an almost complex $n$-dimensional manifold $(M,J)$ and let $E \subset \overline\Omega $  be a smooth totally real $n$-dimensional manifold. Suppose that a function $u$ is plurisubharmonic on $\Omega$ and 
$u^* \vert_{E \setminus \{ p \}} \le C$ for some $p \in E \cap b\Omega$. Then for every $\varepsilon > 0$ there exists a neighborhood $U$ of $p$ such that $u \le C + \varepsilon$ on $A_\alpha(p) \cap U$ for every $\alpha > 1$.
\end{thm}
For the proof it suffices to use properties (b7)-(b9) of the family (\ref{disc2}) and to apply the two constant theorem to the restriction of $u$ on every disc.

\bigskip

\begin{cor}
\label{Lindel2}
Let $\Omega$ be a smoothly bounded domain in an almost complex $n$-dimensional manifold $(M,J)$ and let $E \subset \overline\Omega $  be a smooth totally real $n$-dimensional manifold. Suppose that a function $u$ is plurisubharmonic  and upper bounded on $\Omega$ and 
$\lim\sup_{E \ni q \to p} u^*(q) = -\infty$ for some $p \in E \cap b\Omega$. Then for every  $\alpha > 1$ one has 
$\lim_{A_\alpha(p) \ni q \to p} u(q) = -\infty$.
\end{cor}
Finally, we have

\begin{cor}
\label{Lindel3}
Let $\Omega$ be a smoothly bounded domain in an almost complex $n$-dimensional manifold $(M,J)$ and let $L \in \Omega$ be a generic $(n+1)$-dimensional manifold with the totally real boundary $N \subset b\Omega $. Suppose that a function $u$ is plurisubharmonic and bounded above on $\Omega$ and 
$\lim\sup_{L \ni q \to N} u(q) = -\infty$. Then $u \equiv -\infty$.
\end{cor}
Indeed, for any point $p \in N$ consider a totally real $n$-dimensional manifold $E \subset L \cup \{ p \}$ containing $p$ and apply Corollary \ref{Lindel2}. 

\end{document}